\documentclass{article}
\usepackage{amsmath,amssymb}
\usepackage{amsthm}
\usepackage{graphicx}

\newcounter{exa}
\newcommand{\R}{{\mathbb R}}
\newcommand{\N}{{\mathbb N}}

\newcommand{\F}{{\mathcal F}}
\newcommand{\G}{{\mathcal G}}
\newcommand{\X}{{\mathcal X}}
\newcommand{\E}{\mathbb{E}}
\newcommand{\Var}{\operatorname{Var}}
\newcommand{\cov}{\operatorname{Cov}}

\newcommand{\cqfd}{\hfill $\Box$}

\renewcommand{\L}{\mathcal{L}}

\renewcommand{\H}{{\mathcal H}}

\newtheorem{lemma}{Lemma}
\newtheorem{theorem}{Theorem}

\begin{document}

\title{Empirical processes of iterated maps that contract on average}

\author{Olivier Durieu\\
{\small Universit\'e Fran\c cois-Rabelais de Tours}\\
{\small Laboratoire de Math\'ematiques et Physique Th\'eorique UMR-CNRS 7350,}\\
{\small F\'ed\'eration Denis Poisson FR-CNRS 2964}\\
{\small olivier.durieu@lmpt.univ-tours.fr}}


\date{\today}

\maketitle
\begin{abstract}
We consider a Markov chain obtained by random iterations of Lipschitz maps $T_i$ chosen with a probability $p_i(x)$ depending on the current position $x$.
We assume this system has a property of ``contraction on average'', that is 
\begin{equation*}
 \sum_i d(T_ix,T_iy)p_i(x) < \rho d(x,y)
\end{equation*}
for some $\rho<1$. In the present note, we study the weak convergence of the empirical process associated to this Markov chain.
\end{abstract}



\section{Introduction}

The study of limit behavior of partial sums of random variables has a long time interest in dynamical systems or Markov chains theory. Existence of (attractive) invariant measure, law of large numbers, central limit theorem, almost sure invariance principle are results that characterize the statistical properties of such systems.
In this note, we study the limit behavior of empirical processes, associated to some random iterative Lipschitz models. 
The analyze of empirical processes convergence is also of interest in dynamical systems in order to derive statistical tests, as Kolmogorov-Smirnov test. See for example \cite{ColmarSch04} where expanding maps of the interval are considered. The iterative Lipschitz models that we study are described in Sections \ref{2} and \ref{3}.
Some background on empirical processes is briefly recalled in Section \ref{4}. In Sections \ref{5} and \ref{6}, we show how results of \cite{DehDurTus12} can be applied to get an empirical central limit theorem for the iterative Lipschitz models considered here.

\section{Iterated Lipschitz maps}\label{2}
Let $\mathcal{X}$ be a locally compact metric space, with a countable basis. Denote by $d$ the metric on $\mathcal{X}$.
Let $T_i$, $i\ge 0$, be a sequence of Lipschitz maps from $\mathcal{X}$ to $\mathcal{X}$ and $p_i$, $i\ge 0$, be a sequence of Lipschitz functions from $\mathcal{X}$ to $[0,1]$ such that
for all $x\in\mathcal{X}$, $\sum_{i\ge 0} p_i(x)=1$.
We will consider the Markov chain with state space $\mathcal{X}$ and transition probability $P$ given by
\[
 P(x,A)=\sum_{i\ge 0}p_i(x)1_{A}(T_ix),
\]
for all $x\in\mathcal{X}$ and all Borel subset $A$ of $\mathcal{X}$.
We also denote by $P$ the Markov operator defined for measurable functions $f$ on $\mathcal{X}$ by
\[
 Pf(x)=\int_\mathcal{X}f(y)P(x,dy)=\sum_{i\ge 0}p_i(x)f(T_ix).
\]
We assume the maps $T_i$ contract on average, that is, there exists $\rho\in(0,1)$ such that for all $x,y,z\in\mathcal{X}$,
\begin{equation}\label{ca}
 \sum_{i\ge 0}d(T_ix,T_iy)p_i(z) < \rho d(x,y).
\end{equation}
The case of constant $p_i$ has been studied in \cite{DubFre66}, \cite{BarElt88} (with applications to image encoding) or \cite{HenHer01}. In particular, an empirical CLT has been proved in \cite{WuSha04}.
For variable $p_i$, such systems have been considered in \cite{DoeFor37}, \cite{Kar53} (with applications in learning models), \cite{BarDemEltGer88} (existence of invariant measure), \cite{Pei93} (clt), \cite{Pol01} (Berry-Essen bounds) or \cite{Wal07} (AISP). A lot of concrete examples, that we do not present in this short note, can be found in the articles cited above.

For our general setting, we need the following extra assumptions:
\begin{equation}\label{h0}
\sup_{x,y,z\in\mathcal{X}, y\ne z}\sum_{i\ge 0}\frac{d(T_iy,T_iz)}{d(y,z)}p_i(x)<+\infty.
\end{equation}
\begin{equation}\label{h1}
\sup_{x,y\in\mathcal{X}}\sum_{i\ge0}\frac{d(T_iy,x_0)}{1+d(y,x_0)}p_i(x) < +\infty \mbox{ for some } x_0\in\mathcal{X}. 
\end{equation}
\begin{equation}\label{h2}
\sup_{x\in\mathcal{X}}\sum_{i\ge 0}\frac{d(T_ix,x_0)}{1+d(x,x_0)}\sup_{y,z\in\mathcal{X}, y\ne z}\frac{|p_i(y)-p_i(z)|}{d(y,z)}<+\infty\mbox{ for some } x_0\in\mathcal{X}. 
\end{equation}
Remark that these three conditions are trivially satisfied when the family of maps $T_i$ is finite.

We will also assume that, for any $x, y\in\mathcal{X}$, there exist sequences of integer $(i_n)_{n\ge1}$ and $(j_n)_{n\ge 1}$ such that, for some $x_0\in\mathcal{X}$,
\begin{equation}\label{h4}
 d(T_{i_n}\circ\dots\circ T_{i_1}x,T_{j_n}\circ\dots\circ T_{j_1}y)(1+d(T_{j_n}\circ\dots\circ T_{j_1}x,x_0))\xrightarrow[n\rightarrow+\infty]{} 0,
\end{equation}
with $p_{i_n}(T_{i_{n-1}}\circ\dots\circ T_{i_1}x)\dots p_{i_1}(x)>0$ and  $p_{j_n}(T_{j_{n-1}}\circ\dots\circ T_{j_1}y)\dots p_{i_1}(x)>0$ for all $n\ge 1$.
Remark that this assumption is satisfied when \eqref{ca}--\eqref{h2} hold and the $p_i$ are all strictly positive.

\section{Spaces of Lipschitz functions with weights}\label{3}

We will introduce some Banach spaces, on which the operator $P$ acts with good spectral properties. 
Let $\alpha, \beta$ be non-negative real numbers. We denote by $\H_{_\alpha,\beta}$ the space of continuous functions from $\mathcal{X}$ to $\R$ such that
$\|f\|_{\alpha,\beta}=N_\beta(f)+m_{\alpha,\beta}(f)<+\infty$, where
\[
 N_\beta(f)=\sup_{x\in\mathcal{X}}\frac{|f(x)|}{1+d(x,x_0)^\beta} \;\mbox{ and }\; m_{\alpha,\beta}(f)=\sup_{x,y\in\mathcal{X}, x\ne y}\frac{|f(x)-f(y)|}{d(x,y)^\alpha(1+d(x,x_0)^\beta)} 
\]
for some fixed $x_0\in\mathcal{X}$.
The following statement is easily verifiable.
\begin{lemma}\label{lem}
 If $f,g\in\H_{\alpha,\beta}$ with $\|f\|_\infty=N_0(f)\le 1$ and $\|g\|_\infty\le 1$, then $fg\in\H_{\alpha,\beta}$ and
$\|fg\|_{\alpha,\beta}\le \|f\|_{\alpha,\beta}+\|g\|_{\alpha,\beta}$.
\end{lemma}

The interest of introducing these Banach spaces is the following spectral decomposition of the Markov operator which is due to Peign\'e \cite{Pei93}.
\begin{theorem}[Peign\'e \cite{Pei93}]\label{peigne}
Assume \eqref{ca}--\eqref{h4} hold. Let $\alpha,\beta\in(0,1/2)$ with $\alpha<\beta$. Then 
\begin{enumerate}
\item $P$ operates on $\H_{\alpha,\beta}$.
\item There exists an attractive $P$-invariant probability measure $\nu$ which admits a moment of order one, i.e.\ for all $x_0\in\mathcal{X}$,
$
 \int_\mathcal{X}d(x,x_0)\nu(dx)<\infty.
$
\item There exists a bounded operator $Q$ with spectral radius strictly less than $1$
such that $P=\nu + Q$, with $\nu Q= Q\nu = 0$.
\end{enumerate}
\end{theorem}
As a first consequence, Peign\'e obtained a central limit theorem for the process $(f(X_k))_{k\ge 0}$, where $f$ is any bounded Lipschitz function and $(X_k)_{k\ge 0}$ is the Markov chain starting at a fixed point $x$ with transitions given by $P$.
Another consequence is the almost sure invariance principle which was proved by Walkden \cite{Wal07}. As a result, the functional central limit theorem and the law of iterated logarithm hold.
Here we address the question of a limit theorem for the empirical process associated to the Markov chain starting with distribution $\nu$.

\section{Empirical processes}\label{4}

For a stationary process $(X_k)_{k\ge 0}$ and a class $\F$ of measurable functions from $\mathcal{X}$ to $\R$ which are uniformly bounded, we define the empirical process $(U_n(f))_{f\in\F}\in\ell^\infty(\F)$ by
\[
 U_n(f)=\frac{1}{\sqrt{n}}\sum_{i=0}^{n-1}( f(X_i) - \E(f(X_0))).
\]
The classical cases, which were studied first, correspond to the choice $\mathcal{X}=[0,1]$ and $\F=\{1_{(-\infty,t]} : t\in [0,1]\}$, or more generally, $\mathcal{X}=\R^d$ and $\F=\{1_{(-\infty,t]} : t\in\R^d\}$ for some $d\ge 1$. For these particular cases, results of \cite{DehDur11} can be used to derive a central limit theorem for $(U_n(1_{(-\infty,t]}))_{t\in\R^d}$.
More recently, Dehling, Durieu and Tusche \cite{DehDurTus12} established some results for the abstract case of $\F$-indexed empirical processes. They involve the following definition of bracketing number:

\noindent If $l,u:\X\rightarrow \R$ are two functions satisfying $l(x)\leq u(x)$ for all 
$x\in \X$, the bracket $[l,u]$ is defined by $[l,u]:= \{f:\X\rightarrow \R: l\leq f \leq u \}$.
Given $\varepsilon, A>0, r\ge 1$,  a subset $\G$ of some Banach space with norm $\|\cdot \|$, and a probability law $\mu$, 
we call $[l,u]$ an $(\varepsilon,A,\G,L^r(\mu))$-bracket
if $l,u\in \G$ and
\begin{eqnarray*}
 \|u-l\|_r \leq \varepsilon, & \|u\|\leq A,& \|l\|\leq A,
\end{eqnarray*}
where $\|\cdot \|_r$ denotes the $L^r(\mu)$-norm.
Now, for a class of measurable functions $\F$,
the bracketing number $N(\varepsilon,A,\F,\G,L^r(\mu))$ is defined as 
the smallest number of $(\varepsilon,A,\G, L^r(\mu))$-brackets needed to cover 
$\F$.

\section{Main results}\label{5}

We consider the space $\L=\H_{1,0}$ of bounded Lipschitz functions, equipped with the norm $\|\cdot\|_\L=\|\cdot\|_{1,0}$.
Under an assumption on the bracketing number which allows to approximate functions of the class $\F$ by functions of the spaces $\L$, we can derive the following empirical central limit theorem. 

\begin{theorem}\label{mainthm}
Assume \eqref{ca}--\eqref{h4} hold. Consider the Markov chain $(X_k)_{k\ge0}$ starting with the $P$-invariant distribution $\nu$ and transition $P$.
Let $\F$ be a class of functions from $\mathcal{X}$ to $\R$, uniformly bounded. 
If there exist $\gamma>1$, $C>0$, $r\in(1,2)$, and $c\ge 0$ such that 
\begin{equation*}
 \int_0^1\delta^c \sup_{\varepsilon\ge \delta} N(\varepsilon, \exp(C\varepsilon^{-\frac{1}{\gamma}}),\F,\L,L^r(\nu))d\delta<\infty
\end{equation*}
then the empirical process $(U_n(f))_{f\in\F}$ associated to $(X_k)_{k\ge0}$ converges in distribution in the space $\ell^\infty(\F)$ to a tight centered Gaussian process $(W(f))_{f\in\F}$ with covariances
\begin{align}\label{cov}
 \cov(W(f),W(g))&=\sum_{k=0}^\infty \cov(f(X_0),g(X_k)) + \sum_{k=1}^\infty \cov(f(X_k),g(X_0)).
\end{align}
 
\end{theorem}

This theorem can be applied for several classes of functions $\F$. Examples of such classes, as indicators of balls or ellipsoids, can be found in \cite{DehDurTus12}. Let us state the result corresponding to the classical empirical process indexed by left infinite rectangles in $\R^d$.
\begin{theorem}\label{applithm}
Assume \eqref{ca}--\eqref{h4} hold and $\mathcal{X}=\R^d$. Consider the Markov chain $(X_k)_{k\ge0}$ starting with the $P$-invariant distribution $\nu$ and transition $P$.
If there exist $\gamma>1$ and $C>0$ such that the distribution function $F$ of $\nu$ satisfies
\begin{equation*}\label{mod}
w_F(\delta)\le C |\log(\delta)|^{-\gamma},
\end{equation*}
then the empirical process $(U_n(1_{(-\infty,t]}))_{t\in\R^d}$
converges in distribution in the Skorohod space $D(\R^d)$ to a tight centered Gaussian process $(W(t))_{t\in\R^d}$ with covariances as in \eqref{cov}. 
\end{theorem}

\section{Proofs}\label{6}

\noindent{\em Proof of Theorem \ref{mainthm}.~}
We will apply Theorem 1.1 of \cite{DehDurTus12}.
First, as a consequence of the spectral decomposition of Peign\'e's theorem, we have the central limit theorem for functions of $\L$ (see \cite{Pei93}), i.e.\
for all $f\in\L$, 
\begin{equation}\label{clt}
 \frac{1}{\sqrt{n}}\sum_{k=0}^{n-1}(f(X_k)-\E(f(X_0)))\xrightarrow{\mathcal{D}} \mathcal{N}(0,\sigma^2(f)),
\end{equation}
where $\sigma^2(f)=\Var(f(X_0)) + 2\sum_{k=1}^\infty \cov(f(X_0),f(X_k))$.
Now, set $s=\frac{r}{r-1}\in(2,+\infty)$, $\beta=1/s$ and choose $\alpha<\beta$.
From Theorem \ref{peigne}, we infer that there exist a constant $C>0$ and $\theta\in(0,1)$ such that for all $f\in\H_{\alpha,\beta}$,
\begin{equation}\label{pn}
\|P^nf-\nu(f)\|_{\alpha,\beta}\le C \|f\|_{\alpha,\beta} \theta^n.
\end{equation}
We will deduce a multiple mixing property for $\L$ functions.
Fix a positive integer $p$. 
All along the proof, $C$ always denotes a positive constant which does not depend on the choice of the functions $f$, but its value may change.
Observe that for each $f\in\H_{\alpha,\beta}$,
\begin{align*}
\|f\|_s^s&=\int_{\mathcal{X}} |f(x)|^s \nu(dx)\\
&\le \int_\mathcal{X} \left(\frac{|f(x)|}{1+d(x,x_0)^\beta}\right)^s \left(\frac{1+d(x,x_0)^\beta}{(1+ d(x,x_0))^\beta}\right)^s (1+d(x,x_0)) \nu(dx)\\
&\le C N_\beta(f)^s 
\end{align*}
since $\frac{1+d(x,x_0)^\beta}{(1+ d(x,x_0))^\beta}$ is bounded and $\nu$ has a finite moment of order $1$.
Thus for all $f\in\H_{\alpha,\beta}$, $f\in L^s(\nu)$ and
\begin{equation}\label{ls}
 \|f\|_s\le C\|f\|_{\alpha,\beta}.
\end{equation}
Moreover, for each $f\in\L$ with $\|f\|_\infty\le 1$, we have $f\in\H_{\alpha,\beta}$ and $\|f\|_{\alpha,\beta}\le 3\|f\|_\L$.
For all $i_0<i_1<\dots<i_p$ in $\N$, 
by H\"older's inequality and \eqref{ls}, we have
\begin{align}
&|\cov(f(X_{i_0})\dots f(X_{i_q}),f(X_{i_{q+1}})\dots f(X_{i_p}))|\nonumber\\
&\le \E(|f(X_{i_0})\dots f(X_{i_q})|^r)^{\frac{1}{r}}\nonumber\\
&\qquad \times\E(|\E(f(X_{i_{q+1}})\dots\E(f(X_{i_p})|X_{i_{p-1}})\dots|X_{i_q})-\E(f(X_{i_{q+1}})\dots f(X_{i_p}))|^s)^{\frac{1}{s}} \nonumber\\
&\le \|f\|_r\|P^{i_{q+1}-i_q}g-\nu(g)\|_{\alpha,\beta},\label{eq1}
\end{align}
where $g=fP^{i_{q+2}-i_{q+1}}(f\dots P^{i_p-i_{p-1}}f)$.
Since for all $n$, $P^nf\in\H_{\alpha,\beta}$ and $\|P^nf\|_\infty\le \|f\|_\infty\le 1$, Lemma \ref{lem} shows that
$g\in\H_{\alpha,\beta}$. We can apply \eqref{pn} to get 
\begin{equation}\label{eq2}
\|P^{i_{q+1}-i_q}g-\nu(g)\|_{\alpha,\beta}\le  C\|g\|_{\alpha,\beta}\theta^{i_{q+1}-i_q}.
\end{equation}
Further, the spectral radius of $P$ on $\H_{\alpha,\beta}$ is $1$, thus there exist a constant
$C$ such that $\|P^nh\|_{\alpha,\beta}\le C\|h\|_{\alpha,\beta}$ for all $n$ and all $h\in\H_{\alpha,\beta}$.
Again, by Lemma \ref{lem}, there exists a constant which is independent of $f$ such that
\begin{equation}\label{eq3}
 \|g\|_{\alpha,\beta}\le C\|f\|_{\alpha,\beta}\le C\|f\|_\L.
\end{equation}
Thus, by \eqref{eq1}, \eqref{eq2}, \eqref{eq3}, we get
\[
|\cov(f(X_{i_0})\dots f(X_{i_q}),f(X_{i_{q+1}})\dots f(X_{i_p}))|
\le C\|f\|_r \|f\|_\L\theta^{i_{q+1}-i_q}.
\]
Then $(X_k)_{k\ge 0}$ has a multiple mixing property with exponential decay on the space $\L$. Applying Theorem 4 of \cite{DehDur11}, we obtain the following moment bound: for all
$f\in\L$ with $\|f\|_\infty\le 1$,
\begin{equation}\label{mb}
\E\left(\left(\sum_{k=0}^{n-1}(f(X_k)-\E(f(X_k)))\right)^{2p}\right) \le C\sum_{i=0}^p n^i \|f\|_r^i \log^{2p-i}(\|f\|_\L+1).
\end{equation}
Therefore, thanks to \eqref{clt} and \eqref{mb}, Theorem 1.1 of \cite{DehDurTus12} can be applied. The covariance structure \eqref{cov} of the limit process can be derived from the expression of $\sigma^2(f)$ and \eqref{pn}.
\cqfd

\medskip

\noindent{\em Proof of Theorem \ref{applithm}.~}
By assumption there exists $r\in(1,2)$ such that $\gamma>r$ and \eqref{mod} holds. By Proposition 4.1 of \cite{DehDurTus12}, if $\F=\{1_{(-\infty,t]} : t\in\R^d\}$, we have
\[
 N(\varepsilon,\exp(C\varepsilon^{-\frac{1}{\gamma}}),\F,\L,L^r(\F))=O(\varepsilon^{-dr})
\]
and Theorem \ref{mainthm} completes the proof.
\cqfd

\footnotesize
\bibliographystyle{plain}

\begin{thebibliography}{10}

\bibitem{BarDemEltGer88}
M.~F. Barnsley, S.~G. Demko, J.~H. Elton, and J.~S. Geronimo.
\newblock Invariant measures for {M}arkov processes arising from iterated
  function systems with place-dependent probabilities.
\newblock {\em Ann. Inst. H. Poincar\'e Probab. Statist.}, 24(3):367--394,
  1988.

\bibitem{BarElt88}
Michael~F. Barnsley and John~H. Elton.
\newblock A new class of {M}arkov processes for image encoding.
\newblock {\em Adv. in Appl. Probab.}, 20(1):14--32, 1988.

\bibitem{ColmarSch04}
P.~Collet, S.~Martinez, and B.~Schmitt.
\newblock Asymptotic distribution of tests for expanding maps of the interval.
\newblock {\em Ergodic Theory Dynam. Systems}, 24(3):707--722, 2004.

\bibitem{DehDur11}
Herold Dehling and Olivier Durieu.
\newblock Empirical processes of multidimensional systems with multiple mixing
  properties.
\newblock {\em Stochastic Process. Appl.}, 121(5):1076--1096, 2011.

\bibitem{DehDurTus12}
Herold Dehling, Olivier Durieu, and Marco Tusche.
\newblock Emipirical process of markov chains and dynamical systems indexed by
  classes of functions.
\newblock {\em submitted, arXiv:1201.2256}.

\bibitem{DoeFor37}
Wolfgang Doeblin and Robert Fortet.
\newblock Sur des cha\^{\i}nes \`a liaisons compl\`etes.
\newblock {\em Bull. Soc. Math. France}, 65:132--148, 1937.

\bibitem{DubFre66}
Lester~E. Dubins and David~A. Freedman.
\newblock Invariant probabilities for certain {M}arkov processes.
\newblock {\em Ann. Math. Statist.}, 37:837--848, 1966.

\bibitem{HenHer01}
Hubert Hennion and Lo{\"{\i}}c Herv{\'e}.
\newblock {\em Limit theorems for {M}arkov chains and stochastic properties of
  dynamical systems by quasi-compactness}, volume 1766 of {\em Lecture Notes in
  Mathematics}.
\newblock Springer-Verlag, Berlin, 2001.

\bibitem{Kar53}
Samuel Karlin.
\newblock Some random walks arising in learning models. {I}.
\newblock {\em Pacific J. Math.}, 3:725--756, 1953.

\bibitem{Pei93}
Marc Peign{\'e}.
\newblock Iterated function systems and spectral decomposition of the
  associated {M}arkov operator.
\newblock In {\em Fascicule de probabilit\'es}, volume 1993 of {\em Publ. Inst.
  Rech. Math. Rennes}, page~28. Univ. Rennes I, Rennes, 1993.

\bibitem{Pol01}
M.~Pollicott.
\newblock Contraction in mean and transfer operators.
\newblock {\em Dyn. Syst.}, 16(1):97--106, 2001.

\bibitem{Wal07}
C.~P. Walkden.
\newblock Invariance principles for iterated maps that contract on average.
\newblock {\em Trans. Amer. Math. Soc.}, 359(3):1081--1097 (electronic), 2007.

\bibitem{WuSha04}
Wei~Biao Wu and Xiaofeng Shao.
\newblock Limit theorems for iterated random functions.
\newblock {\em J. Appl. Probab.}, 41(2):425--436, 2004.

\end{thebibliography}

\end{document}